\documentclass[a4paper,10pt]{article}
\usepackage{graphicx}
\usepackage{epsfig}
\usepackage{amsmath,amssymb}
\usepackage[latin1]{inputenc}
\begin{document}
\title{\bf Solution and Asymptotic
Behavior for a Nonlocal Coupled System of Reaction-Diffusion}
\author{
Carlos Alberto Raposo
\thanks{UFSJ, Pra\c{c}a Frei Orlando 170, Cep 36307-352, S\~{a}o
Jo\~{a}o del-Rei, MG, Brasil. raposo@ufsj.edu.br}
\quad
Mauricio Sep\'{u}lveda
\thanks{Departamento de Ingenier\'{\i}a Matem\a'{a}tica,
Universidad de Concepci\a'{o}n, Casilla 160-C, Concepci\a'{o}n,
Chile. mauricio@ing-mat.udec.cl} 
\quad 
Octavio Vera
Villagr\a'{a}n.
\thanks{Departamento de Matem\a'{a}tica, Universidad
del B\a'{\i}o-B\a'{\i}o, Collao 1202, Casilla 5-C, Concepci\a'{o}n,
Chile. overa@ubiobio.cl}
\\
Ducival Carvallo
Pereira\thanks{UFPA, Rua Augusto Corr\^{e}a 01, CEP 66075-110,
Par\'{a}, Brazil. ducival@oi.com.br} 
\quad
Mauro Lima
Santos\thanks{UFPA, Rua Augusto Corr\^{e}a 01, CEP 66075-110,
Par\'{a}, Brazil. ls@ufpa.br}
}

\date{} \maketitle

\begin{abstract}
\noindent This paper concerns with existence, uniqueness and
asymptotic behavior of the solutions for a nonlocal coupled system
of reaction-diffusion. We prove the existence and uniqueness of weak
solutions by the Faedo-Galerkin method and exponential decay of
solutions by the classic energy method. We improve the results
obtained by Chipot-Lovato and Menezes for coupled systems. A
numerical scheme is presented.
\end{abstract}
\renewcommand{\theequation}{\thesection.\arabic{equation}}
\setcounter{equation}{0}\section{Introduction}We consider the
reaction-diffusion coupled system in parallel way via parameter
$\alpha>0$ of the form
\begin{eqnarray}
\label{e101}u_{t} - a(l(u))\,\Delta u + f(u - v) & = & \alpha\,(u -
v)\quad\mbox{in}\quad\Omega\times (0,\,T)
\\
\label{e102}v_{t} - a(l(v))\,\Delta v - f(u - v) & = & \alpha\,(v -
u)\quad \mbox{in}\quad\Omega\times (0,\,T) \\
\label{e103}u = v & = & 0 \quad\mbox{in}\quad \partial\Omega\times
(0,\,T) \\
\label{e104}u(x,\,0) & = &
u_{0}(x)\quad\mbox{in}\quad \Omega \\
\label{e105}v(x,\,0) & = & v_{0}(x)\quad\mbox{in}\quad\Omega
\end{eqnarray}
where $u=u(x,\,t)$ and $v=v(x,\,t)$ are real valued functions.
$\Omega$ is a bounded domain of $\mathbb{R}^{n},$ $\partial \Omega$
is the boundary of $\Omega$ of class $C^{2}.$
$f:\mathbb{R}\rightarrow \mathbb{R}$ and
$a:\mathbb{R}\rightarrow\mathbb{R}$ are Lipschitz continuous
functions with $a(\xi)\geq m > 0.$
$l:L^{2}(\Omega)\rightarrow \mathbb{R}$ is a continuous linear form.\\
For the last several decades, various types of equations have been
employed as some mathematical model describing physical, chemical,
biological and ecological systems. Among them, the most successful
and crucial one is the following model of semilinear parabolic
partial differential equation, called the reaction-diffusion system
\begin{eqnarray}
\label{e106}\frac{\partial u}{\partial t} - D\,\Delta u - f(u)= 0,
\end{eqnarray}
where $f:\mathbb{R}^{n}\rightarrow\mathbb{R}^{n}$ is a nonlinear
function, and $D$ is an $n\times n$ real matrix of diffusion. This
reaction-diffusion model is obtained by combining the system of
nonlinear ordinary differential equations called the reaction system
\begin{eqnarray}
\label{e107}\frac{d u}{d t} - f(u)=0,
\end{eqnarray}
and the system of linear partial differential equation called the
diffusion system
\begin{eqnarray}
\label{e108}\frac{\partial u}{\partial t}- D\,\Delta u = 0.
\end{eqnarray}
In 1998, L. A. F. Oliveira \cite{oliveira} considered the
reaction-diffusion system where $D$ was a $n\times n$ real matrix
and $f:\mathbb{R}^{n} \rightarrow \mathbb{R}^{n}$ is a $C^{2}$
function. He studied the exponential decay for some cases. Except
for some publications on the subject, such as the searching for
traveling waves solutions and some problem in ecology and epidemic
theory, most of authors assume that diffusion matrix $D$ is
diagonal, so that the coupling between the equations  are present
only on the nonlinearity of the reaction term $f.$ However,
cross-diffusion phenomena are not uncommon (see \cite{capasso} and
references therein) can be treated as equations like in which $D$ is
not even diagonalizable. In 1997, M. Chipot and B. Lovat
\cite{chipot} studied the existence and uniqueness of the solutions
for non local problems
\begin{eqnarray}
\label{e109}u_{t} - a\left(\int_{\Omega}u(x,\,t)\,dx\right)\Delta u
& = & f(x,\,t) \quad \mbox{in} \quad
\Omega\times(0,\,T) \\
\label{e110}u(x,\,t)  & = & 0 \quad\mbox{on}\quad\partial\Omega\times(0,\,T)\\
\label{e111}u(x,\,0) & = & u_{0}(x)\quad\mbox{on}\quad\Omega
\end{eqnarray}
where $\Omega$ is a bounded open subset in $\mathbb{R}^{n}$ $n\geq
1$ with smooth boundary $\partial\Omega.$ $T$ is some arbitrary
time. $a$ is some function from $\mathbb{R}$ into $(0,\,+\infty).$
This problem arises in various situations, for instance $u$ could
describe the density of a population(for instance of bacteria)
subject to spreading. The diffusion coefficient $a$ is then supposed
to depend on the entire population in the domain rather than on the
local density i. e. moves are guided by considering the global state
of the medium. They proven the
following result:\\
\\
{\bf Theorem 1.1.} {\it Let $T_{0}>0,$ $u_{0}\in L^{2}(\Omega),$
$u_{0}\geq 0,$ $u_{0}\not\equiv 0.$ Let $a$ be a continuous function
positive in a neighborhood of $\int_{\Omega}u_{0}\,dx.$ Then for
$f\in L^{2}([0,\,T]:\,H^{-1}(\Omega))$ there exists $0<T\leq T_{0}$
and $u$ solution to \eqref{e109}-\eqref{e111} such that
\begin{eqnarray*}
&  & u\in L^{2}([0,\,T]:\,H_{0}^{1}(\Omega))\cap
C^{0}([0,\,T]:\,L^{2}(\Omega)), \\
&  & u_{t}\in L^{2}([0,\,T]:\,H^{-1}(\Omega)),\\
&  & <u_{t},\,v> + \;a\left(\int_{\Omega}u\,dx\right)(\nabla u\cdot
\nabla v) = <f,\,v>\quad \forall\;v\in H_{0}^{1}(\Omega),\quad a.\,
e.\quad t\in [0,\,T_{0}]
\end{eqnarray*}
where $(\nabla u\cdot \nabla v)=\int_{\Omega}\nabla u\cdot\nabla
v\,dx.$}\\
\\
In 2005, S. D. Menezes \cite{Menezes}, give a simple extension of
the result obtained by M. Chipot and B. Lovat \cite{chipot},
considering $a=a(l(u)),$ $f=f(u)$ continuous functions. Indeed, they
studied the existence, uniqueness and periodic solution for the
following parabolic problem
\begin{eqnarray}
\label{e112}u_{t} - a(l(u))\Delta u + f(u) & = & h \quad \mbox{in}
\quad
\Omega\times(0,\,T) \\
\label{e113}u(x,\,t)  & = & 0 \quad\mbox{on}\quad\partial\Omega\times(0,\,T)\\
\label{e114}u(x,\,0) & = & u_{0}(x)\quad\mbox{on}\quad\Omega
\end{eqnarray}
where $\Omega$ is a bounded open subset in $\mathbb{R}^{n},$ $n\geq
1$ with smooth boundary $\partial\Omega.$ $T$ is some arbitrary
time. $l:L^{2}(\Omega)\rightarrow \mathbb{R}$ is a nonlinear form,
$h\in L^{2}(0,\,T:\,H^{-1}(\Omega))$ and $T>0$ is some fixed time.
This problem is nonlocal in the sense that the diffusion coefficient
is determined by a global quantity. This kind of problems, besides
its mathematical motivation because of presence of the nonlocal term
$a(l(u)),$ arises from physical situations related to migration of a
population of bacteria in a container in which the velocity of
migration $\overrightarrow{\nu}=a\,\nabla u$ depends on the global
population in a subdomain $\Omega'\subset\Omega$ given by
$a=a(\int_{\Omega'}u\,dx).$ For more information \cite{chipot} and
reference therein.\\
\\
This article is concerned with to prove the existence, uniqueness
and the exponential decay of the system \eqref{e101}-\eqref{e105}
using the energy method. The method of energy consists of to use
appropriate multipliers to build a functional of Lyapunov, in this
direction we prove that for this types of materials where the
energy, that can flow from one part to another, is strong enough to
produce exponential decay for the solution of the coupled
system.\\
\\
This paper is organized as follows. Before the main result, in
section 2 we briefly outline the notation and terminology to be used
subsequently. In the section three we prove the existence and uniqueness of
solution, in the section four we prove the exponential decay of
solution
of the system.
Finally,  numerical evidence corroborating our theoretical
results is given in section five.  In this paper, we prove the following two theorem:\\
\\
{\bf Theorem 1.2.} {\it Let $(u_{0},\,v_{0})\in L^{2}(\Omega)\times
L^{2}(\Omega)$ and $0<T<+\infty,$ where the time $T$ depends only
$|u_{0}|_{L^{2}(\Omega)}$ and $|v_{0}|_{L^{2}(\Omega)}.$ If
\eqref{e206}-\eqref{e209} holds, then there is at most one solution
of \eqref{e101}-\eqref{e105} in
$L^{2}(0,\,T:\,H_{0}^{1}(\Omega))\cap
C([0,\,T):\,L^{2}(\Omega))\times
L^{2}(0,\,T:\,H_{0}^{1}(\Omega))\cap
C([0,\,T):\,L^{2}(\Omega))$ with initial data $(u(x,\,0),\,v(x,\,0))=(u_{0},\,v_{0}).$}\\
\\
{\bf Theorem 1.3.} {\it Let $(u,\,v)$ be a solution of system
\eqref{e101}-\eqref{e105} given by the theorem 1.2, then there exist
positives constants $C$ and $\eta,$ such that}
\begin{eqnarray}
\label{e115}E(t) \leq  C\,E(0)\,e^{-\eta\,t}.
\end{eqnarray}
\renewcommand{\theequation}{\thesection.\arabic{equation}}
\setcounter{equation}{0}\section{Preliminaries} In this work we
consider the reaction-diffusion coupled system in parallel way via
parameter $\alpha > 0$ as
\begin{eqnarray}
\label{e201}u_{t} - a(l(u))\,\Delta u + f(u - v) & = & \alpha\,(u -
v)\quad\mbox{in}\quad\Omega\times (0,\,T)
\\
\label{e202}v_{t} - a(l(v))\,\Delta v - f(u - v) & = & \alpha\,(v -
u)\quad \mbox{in}\quad\Omega\times (0,\,T) \\
\label{e203}u = v & = & 0 \quad\mbox{in}\quad \partial\Omega\times
(0,\,T) \\
\label{e204}u(x,\,0) & = &
u_{0}(x)\quad\mbox{in}\quad \Omega \\
\label{e205}v(x,\,0) & = & v_{0}(x)\quad\mbox{in}\quad\Omega
\end{eqnarray}
where $\Omega$ is a bounded domain of $\mathbb{R}^{n}$, $\partial
\Omega$ is the boundary of $\Omega$ of class $C^{2}.$
$f:\mathbb{R}\rightarrow \mathbb{R}$ is a Lipschitz continuous
function, that is, there exists $M_{1}>0$ such that
\begin{eqnarray}
\label{e206}|f(s) - f(t)|\leq M_{1}\,|s - t|,\qquad
\forall\;s,\,t\in\mathbb{R}.
\end{eqnarray}
$a:\mathbb{R} \rightarrow\mathbb{R}$ is a Lipschitz continuous
function, that is, there exists $M_{2}>0$ such that
\begin{eqnarray}
\label{e207}|a(s) - a(t)|\leq M_{2}\,|s - t|,\qquad
\forall\;s,\,t\in\mathbb{R}.
\end{eqnarray}
with
\begin{eqnarray}
\label{e208}a(\xi)\geq m > 0,\qquad \forall\;\xi\in\mathbb{R}
\end{eqnarray}
and
\begin{eqnarray}
\label{e209}l:L^{2}(\Omega)\rightarrow \mathbb{R} \quad\mbox{is a
continuous linear form}.
\end{eqnarray}
In the system, the distributed spring coefficient is coupled by the
terms $\alpha\,(u - v)$ and $\alpha\,(v - u).$ In this sense the
Energy can flow from one part to another through
this parameter $\alpha.$\\
\\
By $<\,\cdot\,,\,\cdot\,>$ we will represent the duality pairing
between $X$ and $X',$ $X'$ being the topological dual of the space
$X,$ We represent by $H^{m}(\Omega)$ the usual Sobolev space of
order $m,$ by $H_{0}^{m}(\Omega)$ the closure of
$C_{0}^{\infty}(\Omega)$ in $H^{m}(\Omega),$ and by $L^{2}(\Omega)$
the class of square Lebesgue integrable real functions. In
particular, $H_{0}^{1}(\Omega)$ has inner product
$((\,\cdot\,,\,\cdot\,))$ and norm $||\,\cdot\,||$ given by
\begin{eqnarray*}
((u,\,v)) = \int_{\Omega}\nabla u\cdot \nabla
v\;dx\qquad\mbox{and}\qquad ||u||^{2} = \int_{\Omega}|\nabla
u|^{2}\,dx.
\end{eqnarray*}
For the Hilbert space $L^{2}(\Omega)$ we represent its inner and
norm, respectively, by $(\,\cdot\,,\,\cdot\,)$ and $|\,\cdot\,|,$
defined by
\begin{eqnarray*}
(u,\,v) = \int_{\Omega}u\,v\;dx\qquad\mbox{and}\qquad |u|^{2} =
\int_{\Omega}u^{2}\;dx.
\end{eqnarray*}
Throughout this paper  $c$ is a generic constant, not necessarily
the same at each occasion(it will change from line to line), which
depends in an increasing way on the indicated quantities.\\
\\
We take the initial conditions as following
\begin{eqnarray}
\label{e210}(u_{0}(x),\,v_{0}(x))\in L^{2}(\Omega)\times
L^{2}(\Omega).
\end{eqnarray}
We denote the potential energy associated to this system  by
\begin{eqnarray}
\label{e211}E(t)= \frac{1}{2}\int_{\Omega}\left[\,|u|^{2} + |v|^{2}
\,\right]\,dx.
\end{eqnarray}
\renewcommand{\theequation}{\thesection.\arabic{equation}}
\setcounter{equation}{0}\section{Existence and Uniqueness of a local
solution} In this section, we will prove that for
$(u_{0},\,v_{0})\in L^{2}(\Omega)\times L^{2}(\Omega)$ there exists
a unique solution of \eqref{e206}-\eqref{e209} in
$L^{2}(0,\,T:\,H_{0}^{1}(\Omega))\cap
C([0,\,T):\,L^{2}(\Omega))\times
L^{2}(0,\,T:\,H_{0}^{1}(\Omega))\cap C([0,\,T):\,L^{2}(\Omega))$
where the time $T$ depends only $|u_{0}|_{L^{2}(\Omega)}$ and
$|v_{0}|_{L^{2}(\Omega)}.$ We make use of Faedo-Galerkin
approximation for to prove the existence of weakly solutions. We
write the system \eqref{e201}-\eqref{e205} in the following form,
\begin{eqnarray}
\label{e301}u_{t} - a(l(u))\,\Delta u  & = & -\;f(u,\,v) + g(u,\,v)
\quad\mbox{in}\quad Q = \Omega\times(0,\,T)\\
\label{e302}v_{t} - a(l(v))\,\Delta v  & = &  f(u,\,v) - g(u,\,v)
\quad\mbox{in}\quad Q = \Omega\times(0,\,T)\\
\label{e303}u = v & = & 0 \quad\mbox{on}\quad\partial\Omega\times (0,\,T)\\
\label{e304}u(x,\,0) & = & u_{0}(x)\quad \mbox{in}\quad\Omega\\
\label{e305}v(x,\,0) & = & v_{0}(x)\quad \mbox{in}\quad\Omega
\end{eqnarray}
where we denote $f(u,\,v)\equiv f(u - v)$ and $g(u,\,v)\equiv \alpha\,(u - v).$\\
\\
{\bf Theorem 3.1}(Existence). {\it Let $(u_{0},\,v_{0})\in
L^{2}(\Omega)\times L^{2}(\Omega)$ and $0<T<+\infty.$ If
\eqref{e206}-\eqref{e209} holds, then there exists $(u,\,v)$
solution of \eqref{e301}-\eqref{e305} such that }
\begin{eqnarray}
\label{e306}&  & (u,\,v)\in L^{2}(0,\,T:\,H_{0}^{1}(\Omega))\cap
C([0,\,T):\,L^{2}(\Omega))\times
L^{2}(0,\,T:\,H_{0}^{1}(\Omega))\cap
C([0,\,T):\,L^{2}(\Omega))\\
\label{e307}&  & (u_{t},\,v_{t})\in
L^{2}(0,\,T:\,H^{-1}(\Omega))\times L^{2}(0,\,T:\,H^{-1}(\Omega))
\end{eqnarray}
\begin{eqnarray}
\label{e308}\frac{d}{dt}(u,\,h_{1}) + a(l(u))\int_{\Omega}\nabla
u\cdot\nabla h_{1}\,dx = -\int_{\Omega}f(u,\,v)\,h_{1}\,dx +
\int_{\Omega}g(u,\,v)\,h_{1}\,dx
\end{eqnarray}
{\it for all $h_{1}\in H_{0}^{1}(\Omega),$ where \eqref{e308} must
be understood as an equality in $\mathcal{D}'(0,\,T).$}
\begin{eqnarray}
\label{e309}\frac{d}{dt}(v,\,h_{2}) + a(l(v))\int_{\Omega}\nabla
v\cdot\nabla h_{2}\, dx = \int_{\Omega}f(u,\,v)\,h_{2}\,dx -
\int_{\Omega}g(u,\,v)\,h_{2}\,dx
\end{eqnarray}
{\it for all $h_{2}\in H_{0}^{1}(\Omega),$ where \eqref{e309} must
be understood as an equality in $\mathcal{D}'(0,\,T).$} \\
\\
{\it Proof.}\quad i) {\bf Approximate problem:} Let
$\{w_{j}\}_{j\in\mathbb{N}}$ be a Hilbertian basis of
$H_{0}^{1}(\Omega)$ (cf. H. Brezis, \cite{br1}). For each
$j=1,\,2,\,3,\,\ldots$ represent by $\mathbb{V}_{j},$ the subspace
of $H_{0}^{1}(\Omega)$ generated by
$\{w_{1}\,,w_{2},\,\ldots,\,w_{j}\}.$ The approximate problem,
associated with \eqref{e301}-\eqref{e305}, consists of to find
$u_{j},\,v_{j}\in \mathbb{V}_{j}$ such that
\begin{eqnarray}
\label{e310}&  & (u_{j}',\,h_{1}) - a(l(u_{j}))(\Delta
u_{j},\,h_{1}) = -\;
(f(u_{j},\,v_{j}),\,h_{1})+(g(u_{j},\,v_{j}),\,h_{1}),
\quad\forall\;h_{1}\in \mathbb{V}_{j}\\
\label{e311}&  & (v_{j}',\,h_{2}) - a(l(v_{j}))(\Delta
v_{j},\,h_{2}) =
(f(u_{j},\,v_{j}),\,h_{2})-(g(u_{j},\,v_{j}),\,h_{2}),
\quad\forall\;h_{2}\in \mathbb{V}_{j}\\
\label{e312}&  & u_{j}(0)=u_{0j}\rightarrow u_{0}\quad\mbox{strongly
in}\quad L^{2}(\Omega) \\
\label{e313}&  & v_{j}(0)=v_{0j}\rightarrow v_{0}\quad\mbox{strongly
in}\quad L^{2}(\Omega)
\end{eqnarray}
Let $h_{1}=w_{i}(x)$ and $h_{2}=w_{i}(x)$ for $i=1,\,\ldots,\,j,$
then in \eqref{e310}-\eqref{e313} we have for
$\theta_{kj},\,\phi_{kj}\in C^{\infty}(\Omega)$
\begin{eqnarray*}
&  & \left(\sum_{k=1}^{j}\theta_{kj}'(t)w_{k}(x),\,w_{i}(x)\right) -
a \left(l\left(\sum_{k=1}^{j}\theta_{kj}(t)w_{k}(x)\right)\right)
\left(\Delta\sum_{k=1}^{j}
\theta_{kj}(t)w_{k}(x),\,w_{i}(x)\right)\\
& = & -\left(f\left(\sum_{k=1}^{j}\theta_{kj}(t)w_{k}(x),\,
\sum_{k=1}^{j}\phi_{kj}(t)w_{k}(x)\right),\,w_{i}(x)\right) \\
&  & +\left(g\left(\sum_{k=1}^{j}\theta_{kj}(t)w_{k}(x),\,
\sum_{k=1}^{j}\phi_{kj}(t)w_{k}(x)\right),\,w_{i}(x)\right)
\end{eqnarray*}
and
\begin{eqnarray*}
&  & \left(\sum_{k=1}^{j}\phi_{kj}'(t)w_{k}(x),\,w_{i}(x)\right) -
a\left(l\left(\sum_{k=1}^{j}\phi_{kj}(t)w_{k}(x)\right)\right)
\left(\Delta\sum_{k=1}^{j}
\phi_{kj}(t)w_{k}(x),\,w_{i}(x)\right) \\
& = & \left(f\left(\sum_{k=1}^{j}\theta_{kj}(t)w_{k}(x),\,
\sum_{k=1}^{j}\phi_{kj}(t)w_{k}(x)\right),\,w_{i}(x)\right) \\
&  & -\left(g\left(\sum_{k=1}^{j}\theta_{kj}(t)w_{k}(x),\,
\sum_{k=1}^{j}\phi_{kj}(t)w_{k}(x)\right),\,w_{i}(x)\right)
\end{eqnarray*}
that is,
\begin{eqnarray}
\label{e314}\theta_{kj}'(t) -
\lambda_{k}\,a(l(u_{j}))\,\theta_{kj}(t) & = &
-\;(f(u_{j},\,v_{j}),\,w_{i}) + (g(u_{j},\,v_{j}),\,w_{i})\\
\label{e315}\phi_{kj}'(t) - \lambda_{k}\,a(l(v_{j}))\,\phi_{kj}(t) &
= & (f(u_{j},\,v_{j}),\,w_{i}) - (g(u_{j},\,v_{j}),\,w_{i}).
\end{eqnarray}
ii) {\bf Approximate solutions:} We will just work with the equation
\eqref{e314}. For \eqref{e315} the result is similar. For
$i,\,k=1,\,\ldots,\,j$ in \eqref{e314}, we have the following system
\begin{displaymath}
\left[
  \begin{array}{c}
    \theta_{1j}' \\
    \theta_{2j}' \\
    \vdots \\
    \theta_{jj}'
  \end{array}
\right]=\left[
  \begin{array}{ccc}
    \lambda_{1}\,a(l(u_{j})) &  & 0 \\

          & \ddots&  \\
    0 &  & \lambda_{j}\,a(l(u_{j})) \\
  \end{array}
\right]\left[
          \begin{array}{c}
            \theta_{1j} \\
            \theta_{2j} \\
            \vdots \\
            \theta_{jj}\\
          \end{array}
        \right]
        -\left[\begin{array}{c}
                    (f,\,w_{1}) \\
                    (f,\,w_{2}) \\
                    \vdots \\
                    (f,\,w_{j}) \\
                  \end{array}
                \right]
        +\left[
                  \begin{array}{c}
                    (g,\,w_{1}) \\
                    (g,\,w_{2}) \\
                    \vdots \\
                    (g,\,w_{j}) \\
                  \end{array}
                \right]
\end{displaymath}
that is,
\begin{eqnarray}
\label{e316}X' & = & F(X,\,t)\\
\label{e317}X(0) & = & X_{0}
\end{eqnarray}
where $F(X,\,t)=A\,X + B$ and $X_{0}=[
\alpha_{1j},\,\alpha_{2j},\,\ldots,\,\alpha_{jj}]^{T}.$ The system
\eqref{e316}-\eqref{e317} is equivalent to system of ordinary
differential equations of first order. Let us show that the system
\eqref{e316}-\eqref{e317} is in
the conditions of Carath\'eodory's theorem.\\
\\
{\it Claim.} For fixed  $X$, we will show that $A$ and $B$ are
measurable in $t.$\\
\\
In fact, we observed that the matrix  $A$ is formed for the elements
$\lambda_{k}\,a(l(u_{j}))$ with $k=1,\,2,\,\ldots,\,j.$ Since $l$ is
a lineal and continuous form and the operator $a$ is continuous,
then the composition $a(l(u_{j})) $ is also continuous; therefore
$\lambda_{k}\,a(l(u_{j}))$  is continuous for $k=1,\,2,\,\ldots,\,j$
and then $A$ is measurable in $t.$ On the other hand, let us observe
that $B $ is formed by the elements $(f(u_{j},\,v_{j}),\,w_{i})$ and
$(g(u_{j},\,v_{j}),\,w_{i}),$ with $i=1,\,2,\,\ldots,\,j.$ Since $f$
and $g$ are continuous and
$w_{i}\in H_{0}^{1}(\Omega),$ we concludes that $B$ is continuous and
therefore measurable.\\
\\
{\it Claim.} For fixed $t,$ we will show that $F$ is continuous  in
$X.$\\
\\
In fact, notice that $B$ is continuous in $X,$ because $B$ is
constant in relation to $X.$ For continuity of $A\,X,$ is enough we
verify that $A$ is continuous in $X.$ Let $\prod_{k}X =
\theta_{kj}\,(k=1,\,2,\,\ldots,\,j)$ be the projection
$\mathbb{R}^{j}\longrightarrow \mathbb{R}$ and
$\sigma(X)=\prod_{k}X\,w_{k}.$ For each $t$ fixed, as
\begin{eqnarray*}
u_{j}(t)=\sum_{k=1}^{j}\theta_{kj}(t)\,w_{k},
\end{eqnarray*}
we can consider the function
\begin{eqnarray*}
X\longrightarrow a(l(u_{j}))=
a\left(l\left(\sum_{k=1}^{j}\theta_{kj}(t)\,w_{k}\right)\right)=
a\left(l\left(\sum_{k=1}^{j}\prod_{k}X\,w_{k}\right)\right).
\end{eqnarray*}
Since $A$ is lineal combination of continuous functions, proceeds
that $A$ is continuous in $X,$ hence the function $F(X,\,t)$ is
continuous
in $X.$\\
\\
{\it Claim.} Let $K$ be a compact of $\mathbb{D}=\mathbb{E}\times
[0,\,T],$ where $\mathbb{E}=\{X\in \mathbb{R}^{j\times1}:\;
||X||_{\mathbb{R}^{j\times1}}\leq \delta,\;\delta > 0\}.$ We will
show that exists a real function $m_{r}(t),$ integrable in
$[0,\,T],$ so that
\begin{eqnarray*}
||F(X,\,t)||_{\mathbb{R}^{j\times1}}\leq m_{r}(t),\qquad\forall\;
(X,\,t)\in D.
\end{eqnarray*}
In fact, we denote by $||\,\cdot\,||_{p\,q}$ the norm of maximum in
$\mathbb{R}^{pq}$. But $F(X,\,t)=A\,X + B,$ then
\begin{eqnarray*}
||F(X,\,t)||_{j\times1}\leq ||A||_{j\times j}\;||X||_{j\times1} +
||B||_{j\times1}.
\end{eqnarray*}
Since $X\in \mathbb{E},$ we have $||X||_{j\times1}\leq\delta.$ Then
\begin{eqnarray*}
||F(X,\,t)||_{j\times1}\leq \delta\;||A||_{j\times j} +
||B||_{j\times1}.
\end{eqnarray*}
Notice that $\lambda_{k}\,a(l(u_{j}))$  are continuous
functions, hence $||A||_{j\times j}\leq C\,\quad(C>0).$\\
\\
On the other hand, for the matrix $B$ we have
\begin{eqnarray*}
|(f(u_{j},\,v_{j}),\,w_{i})|\leq
|f(u_{j},\,v_{j})|\;|w_{i}| = |f(u_{j},\,v_{j})|,\\
|(g(u_{j},\,v_{j}),\,w_{i})|\leq |g(u_{j},\,v_{j})|\;|w_{i}| =
|g(u_{j},\,v_{j})|.
\end{eqnarray*}
Thus $||F(X,\,t)||_{j\times1}\leq \delta\;C + |f(u_{j},\,v_{j})| +
|g(u_{j},\,v_{j})|\equiv
m_{r}(t),$ where $m_{r}(t)$ is integrable in $[0,\,T].$\\
\\
Hence, the system \eqref{e316}-\eqref{e317} satisfies the conditions
of Carath\'eodory, and then  exists $\{u_{j}(t),\,v_{j}(t)\}\in
[0,\,t_{j})\times[0,\,t_{j}),$ $t_{j}< T_{0}.$ \\
\\
We now have to establish an estimate that permits to extend the
solution $\{u_{j}(t),\,v_{j}(t)\}$ to the whole interval $[0,\,T].$\\
\\
From now on, $\{C_{i}\}_{i=1\ldots7},$ will denote
positive constants, independents of $j$ and $t.$\\
\\
iii) {\bf A priori estimates:} We put $h_{1}=u_{j}$ and
$h_{2}=v_{j}$ in the equations \eqref{e310} and \eqref{e311}
respectively, we have
\begin{eqnarray}
\label{e318}(u_{j}',\,u_{j}) - a(l(u_{j}))\,(\Delta u_{j},\,u_{j})=
-\;
(f(u_{j},\,v_{j}),\,u_{j}) + (g(u_{j},\,v_{j}),\,u_{j})\\
\label{e319}(v_{j}',\,v_{j}) - a(l(v_{j}))\,(\Delta v_{j},\,v_{j})=
(f(u_{j},\,v_{j}),\,u_{j}) - (g(u_{j},\,v_{j}),\,u_{j}).
\end{eqnarray}
Using the boundary condition and the first Green's identity we have
\begin{eqnarray*}
(-\Delta u_{j},\,u_{j})=\int_{\Omega}(-\Delta u_{j})\,u_{j}\,dx =
\int_{\Omega}\nabla u_{j}\cdot\nabla u_{j}\,dx=|\nabla
u_{j}|^{2}=||u_{j}||^{2}.
\end{eqnarray*}
Then, we can write \eqref{e318} as
\begin{eqnarray}
\label{e320}\frac{1}{2}\frac{d}{dt}|u_{j}|^{2} +
a(l(u_{j}))\,||u_{j}||^{2}= -\;(f(u_{j},\,v_{j}),\,u_{j}) +
(g(u_{j},\,v_{j}),u_{j}).
\end{eqnarray}
In a similar way we can write \eqref{e319} as
\begin{eqnarray}
\label{e321}\frac{1}{2}\frac{d}{dt}|v_{j}|^{2} +
a(l(v_{j}))\,||v_{j}||^{2}= (f(u_{j},\,v_{j}),\,u_{j}) -
(g(u_{j},v_{j}),\,u_{j}).
\end{eqnarray}
Adding the equations \eqref{e320} with \eqref{e321} and integrating
of $0$ to $t,$ we obtain
\begin{eqnarray*}
\int_{0}^{t}\frac{1}{2}\frac{d}{ds}\left[\,|u_{j}|^{2} +
|v_{j}|^{2}\,\right]ds +
\int_{0}^{t}\left[\,a(l(u_{j}))\,||u_{j}||^{2} +
a(l(v_{j}))\,||v_{j}||^{2}\,\right]ds = 0,
\end{eqnarray*}
that is,
\begin{eqnarray*}
|u_{j}(t)|^{2} + |v_{j}(t)|^{2} +
2\int_{0}^{t}\left[\,a(l(u_{j}))\,||u_{j}||^{2}+
a(l(v_{j}))\,||v_{j}||^{2}\,\right]ds = |u_{j}(0)|+ |v_{j}(0)|.
\end{eqnarray*}
In the last identity, using \eqref{e206} and \eqref{e208} we obtain
\begin{eqnarray}
\label{e322}|u_{j}(t)|^{2} + |v_{j}(t)|^{2}+
2m\int_{0}^{t}(\,||u_{j}||^{2} + ||v_{j}||^{2}\,)\,ds \leq
|u_{j}(0)|^{2} + |v_{j}(0)|^{2}.
\end{eqnarray}
Since $u_{j}(0)\rightarrow u_{0}$ and $v_{j}(0)\rightarrow v_{0}$
strongly in $L^{2}(\Omega)$ follows that $|u_{j}(0)|^{2} +
|v_{j}(0)|^{2}\leq C.$ Hence
\begin{eqnarray*}
|u_{j}(t)|^{2} + |v_{j}(t)|^{2}\leq C.
\end{eqnarray*}
From where follows that $u_{j}(t)$ and $v_{j}(t)$ are bounded in
$L^{\infty}(0,\,T:\,L^{2}(\Omega)).$ Thus,
\begin{eqnarray*}
\int_{0}^{t}(||u_{j}||^{2} + ||v_{j}||^{2})\,ds\leq C,
\end{eqnarray*}
then $u_{j}(t)$ and $v_{j}(t)$ are limited in $L^{2}(0,\,T:\,
H^{1}_{0}(\Omega)).$\\
\\
From \eqref{e310}-\eqref{e311}, we have that
\begin{eqnarray*}
u_{j}' = a(l(u_{j}))\Delta u_{j} -
f(u_{j},\,v_{j}) + g(u_{j},\,v_{j}) \in H^{-1}(\Omega) \\
v_{j}' = a(l(v_{j}))\Delta v_{j} + f(u_{j},\,v_{j}) -
g(u_{j},\,v_{j})\in H^{-1}(\Omega).
\end{eqnarray*}
Notice that $-\,a(l(u_{j}))\,\Delta u_{j}$ defines an element of
$H^{-1}(\Omega),$ given by the duality
\begin{eqnarray*}
\langle -\,a(l(u_{j}))\,\Delta u_{j},\,h_{1}\rangle =
a(l(u_{j}))\int_{\Omega}\nabla u_{j}\cdot\nabla
h_{1}\,dx,\quad\forall\;h_{1}\in H_{0}^{1}(\Omega).
\end{eqnarray*}
In a similar way we have
\begin{eqnarray*}
\langle -\,a(l(v_{j}))\,\Delta v_{j},\,h_{2}\rangle =
a(l(v_{j}))\int_{\Omega}\nabla v_{j}\cdot\nabla h_{2}\,dx,
\quad\forall\;h_{2}\in H_{0}^{1}(\Omega).
\end{eqnarray*}
Using the fact that $-\,a(l(u_{j}))\,\Delta u_{j},$
$-\,a(l(v_{j}))\,\Delta v_{j}\in H^{-1}(\Omega),$ the dual of
$H_{0}^{1}(\Omega),$ then
they are lineal and continuous forms and therefore both are bounded.\\
\\
Since $u_{j},\,v_{j}\in L^{2}(0,\,T:\,L^{2}(\Omega)),$ then
\begin{eqnarray*}
\int_{\Omega}|f(u_{j},\,v_{j})|\,dx\leq \int_{\Omega}\beta\,|u_{j} -
v_{j}|\,dx \leq \int_{\Omega}\beta\,(|u_{j}| + |v_{j}|)\,dx\\
\leq C\left[\left(\int_{\Omega}|u_{j}|^{2}\,dx\right)^{1/2} +
\left(\int_{\Omega}|v_{j}|^{2}\,dx\right)^{1/2}\right]
\end{eqnarray*}
and
\begin{eqnarray*}
\int_{\Omega}|g(u_{j},\,v_{j})|\,dx = \int_{\Omega}\alpha\,|u_{j}
- v_{j}|\,dx \leq \int_{\Omega}\alpha\,(\,|u_{j}| + |v_{j}|\,)\,dx \\
\leq C\left[\left(\int_{\Omega}|u_{j}|^{2}\,dx\right)^{1/2} +
\left(\int_{\Omega}|v_{j}|^{2}\,dx\right)^{1/2}\right].
\end{eqnarray*}
Therefore
\begin{eqnarray*}
f(u_{j},\,v_{j}),\;g(u_{j},\,v_{j})\in L^{2}(0,\,T:\,
L^{2}(\Omega))\hookrightarrow L^{1}(0,\,T:\,L^{2}(\Omega))
\end{eqnarray*}
and we concludes that  $u_{j}',\;v_{j}'$ are bounded in
$L^{2}(0,\,T:\,H^{-1}(\Omega)).$\\
\\
iv) {\bf Passage to the limit:} We have that
\begin{eqnarray}
\label{e323}&  & u_{j},\;v_{j}\qquad\mbox{are bounded in}\quad
L^{\infty}(0,\,T:\,L^{2}(\Omega))\cap L^{2}(0,\,T:\,H^{1}_{0}(\Omega)),\\
\label{e324}&  & u_{j}',\;v_{j}'\qquad\mbox{are bounded in}\quad
L^{2}(0,\,T:\,H^{-1}(\Omega)).
\end{eqnarray}
Now, due the corollary of Banach-Alouglu (See \cite{ru1}, p. 66), we
can extract subsequences of $u_{j_{k}}\equiv u_{j}$ and
$v_{j_{k}}\equiv v_{j}$ (which we denote with the same symbol) so
that
\begin{eqnarray}
\label{e325}u_{j}\stackrel{\star}{\rightharpoonup}u \quad\mbox{weak
star in}\quad
L^{\infty}(0,\,T:\,L^{2}(\Omega)) \\
\label{e326}v_{j}\stackrel{\star}{\rightharpoonup}v \quad\mbox{weak
star in}\quad
L^{\infty}(0,\,T:\,L^{2}(\Omega)) \\
\label{e327}u_{j}\rightharpoonup u \quad\mbox{weak
 in}\quad L^{2}(0,\,T:\,H^{1}_{0}(\Omega)) \\
\label{e328}v_{j}\rightharpoonup v \quad\mbox{weak
 in}\quad L^{2}(0,\,T:\,H^{1}_{0}(\Omega)).
\end{eqnarray}
Consequently
\begin{eqnarray}
\label{e329}\int_{0}^{T}(u_{j},\,h_{1})\,dt\longrightarrow
\int_{0}^{T}(u,\,h_{1})\,dt,
\qquad\forall \;h_{1}\in L^{\infty}(0,\,T:\,L^{2}(\Omega)) \\
\label{e330}\int_{0}^{T}(u_{j},\,h_{1})\,dt \longrightarrow
\int_{0}^{T}(u,\,h_{1})\,dt, \qquad\forall\; h_{1}\in
L^{2}(0,\,T:\,H^{1}_{0}(\Omega)) \\
\label{e331}\int_{0}^{T}(v_{j},\,h_{2})\,dt \longrightarrow
\int_{0}^{T}(v,\,h_{2})\,dt,
\qquad\forall\; h_{2}\in L^{\infty}(0,\,T:\,L^{2}(\Omega)) \\
\label{e332}\int_{0}^{T}(v_{j},\,h_{2})\,dt\longrightarrow
\int_{0}^{T}(v,\,h_{2})\,dt, \qquad\forall\; h_{2}\in
L^{2}(0,\,T:\,H^{1}_{0}(\Omega))
\end{eqnarray}
For \eqref{e324} it proceeds
\begin{eqnarray}
\label{e333}u_{j}'\rightharpoonup u'\qquad\mbox{weakly in}\quad
L^{2}(0,\,T:\,
H^{-1}(\Omega)) \\
\label{e334}v_{j}'\rightharpoonup v'\qquad\mbox{weakly in}\quad
L^{2}(0,\,T:\,H^{-1}(\Omega)).
\end{eqnarray}
On the other hand, $H_{0}^{1}(\Omega)\stackrel{c}\hookrightarrow
L^{2}(\Omega)\hookrightarrow H^{-1}(\Omega).$ By Lions-Aubin's
compactness Theorem \cite{li1} follows that
\begin{eqnarray}
\label{e335}u_{j}\longrightarrow u \quad\mbox{strongly in}\quad
L^{2}(0,\,T:\,
L^{2}(\Omega)) \\
\label{e336}v_{j}\longrightarrow v \quad\mbox{strongly in}\quad
L^{2}(0,\,T:\,L^{2}(\Omega)).
\end{eqnarray}
The convergence \eqref{e325}-\eqref{e326} means that
\begin{eqnarray*}
\int_{0}^{T}(u_{j}(t),\,w(t))\,dt\longrightarrow
\int_{0}^{T}(u(t),\,w(t))\,dt,\qquad \forall\;w\in
L^{1}(0,\,T:\;L^{2}(\Omega)) \\
\int_{0}^{T}(v_{j}(t),\,w(t))\,dt\longrightarrow
\int_{0}^{T}(v(t),\,w(t))\,dt,\qquad \forall\;w\in
L^{1}(0,\,T:\;L^{2}(\Omega)).
\end{eqnarray*}
We choose $w=\theta\,h_{1}$ with $\theta\in \mathcal{D}(0,\,T),$
$h_{1}\in L^{2}(\Omega)$ and we will show that for all $\theta\in
\mathcal{D}(0,\,T)$ and for all $h_{1}\in L^{2}(\Omega),$
\begin{eqnarray*}
\int_{0}^{T}[\,(g(u_{j},\,v_{j}),\,h_{1}) -
(g(u,\,v),\,h_{1})\,]\,\theta(t)\,dt\longrightarrow 0.
\end{eqnarray*}
Let $T$ be a positive number such that $supp(\theta)\subset
[0,\,T],$ then
\begin{eqnarray*}
\int_{0}^{T}[\,(g(u_{j},\,v_{j}),\,h_{1}) -
(g(u,\,v),\,h_{1})\,]\,\theta(t)\,dt = \int_{0}^{T}(g(u_{j},\,v_{j})
- g(u,\,v),\,h_{1})\,\theta(t)\,dt.
\end{eqnarray*}
Hence, by straightforward calculations
\begin{eqnarray*}
\int_{0}^{T}(g(u_{j},\,v_{j}) - g(u,\,v),\,h_{1})\,\theta(t)\,dt &
\leq &
\int_{0}^{T}\int_{\Omega}|g(u_{j},\,v_{j}) - g(u,\,v)|\;|h_{1}|\;|\theta(t)|\,dx\,dt \\
& = &
\int_{0}^{T}\int_{\Omega}|\alpha\,(u_{j} - v_{j})- \alpha\,(u - v)|\;|h_{1}|\;|\theta(t)|\,dx\,dt \\
& = &
\int_{0}^{T}\int_{\Omega}|\alpha\,(u_{j} - u) - \alpha\,(v_{j} - v)|\;|h_{1}|\;|\theta(t)|\,dx\,dt \\
& \leq & \int_{0}^{T}\int_{\Omega}
(\,\alpha\,|u_{j} - u| + \alpha\,|v_{j} - v|\,)\,|h_{1}|\;|\theta(t)|\,dx\,dt. \\
\end{eqnarray*}
Using $L^{2}(0,\,T:\,L^{2}(\Omega))\hookrightarrow
L^{1}(0,\,T:\,L^{2}(\Omega)) $ and the Cauchy-Schwartz inequality we
obtain
\begin{eqnarray*}
\int_{0}^{T}(g(u_{j},\,v_{j}) - g(u,\,v),\,h_{1})\,\theta(t)\,dt &
\leq & C\int_{0}^{T}\left(\int_{\Omega}|u_{j} -
u|^{2}dx\right)^{1/2}
\left(\int_{\Omega}|h_{1}|^{2}\,dx\right)^{1/2}dt  \\
&  & +\;C\int_{0}^{T}\left(\int_{\Omega}|v_{j} -
v|^{2}dx\right)^{1/2}
\left(\int_{\Omega}|h_{1}|^{2}dx\right)^{1/2}dt.
\end{eqnarray*}
Applying the Cauchy-Schwartz inequality and considering the
convergence \eqref{e335} we obtain
\begin{eqnarray*}
\lefteqn{C\int_{0}^{T}\left(\int_{\Omega}|u_{j} -
u|^{2}dx\right)^{1/2}
\left(\int_{\Omega}|h_{1}|^{2}dx\right)^{1/2}dt } \\
& \leq & C\left(\int_{0}^{T}\int_{\Omega}|u_{j} -
u|^{2}dxdt\right)^{1/2}
\left(\int_{0}^{T}\int_{\Omega}|h_{1}|^{2}dxdt\right)^{1/2}<
\varepsilon.
\end{eqnarray*}
In a similar way using the convergence \eqref{e336} we have
\begin{eqnarray*}
\lefteqn{C\int_{0}^{T}\left(\int_{\Omega}|v_{j} -
v|^{2}dx\right)^{1/2}
\left(\int_{\Omega}|h_{1}|^{2}dx\right)^{1/2}dt } \\
& \leq & C\left(\int_{0}^{T}\int_{\Omega}|v_{j} -
v|^{2}dx\,dt\right)^{1/2}
\left(\int_{0}^{T}\int_{\Omega}|h_{1}|^{2}dx\,dt\right)^{1/2}<
\varepsilon.
\end{eqnarray*}
Therefore we have
\begin{eqnarray*}
\int_{0}^{T}(g(u_{j},\,v_{j}) - g(u,\,v),\,h_{1})\,\theta(t)\,dt <
\varepsilon.
\end{eqnarray*}
Performing similar calculations we can to prove that
\begin{eqnarray*}
\int_{0}^{T}(f(u_{j},\,v_{j}) - f(u,\,v),\,h_{1})\,\theta(t)\,dt <
\varepsilon.
\end{eqnarray*}
We will show now, that for every $\theta\in
\mathcal{D}\left([0,\,T]\right)$ and for every $h_{1}\in
L^{2}(\Omega)$
\begin{eqnarray}
\label{e337}a(l(u_{j}))\int_{0}^{T}\int_{\Omega}\nabla
u_{j}\cdot\nabla h_{1}\;\theta(t)\,dt\longrightarrow
a(l(u))\int_{0}^{T}\int_{\Omega}\nabla u\cdot\nabla
h_{1}\;\theta(t)\,dt.
\end{eqnarray}
It is enough we prove that
\begin{eqnarray}
\label{e338}a(l(u_{j}))\longrightarrow a(l(u))\quad\mbox{in}\quad
L^{2}(0,\,T),\qquad\forall\;T>0.
\end{eqnarray}
Since $a$ is continuous, we will show that
\begin{eqnarray}
\label{e339}l(u_{j})\longrightarrow l(u)\quad\mbox{strongly in}\quad
L^{2}(0,\,T).
\end{eqnarray}
In fact, because
\begin{eqnarray*}
\int_{0}^{T}|l(u_{j}) - l(u))|^{2}dt = \int_{0}^{T}|l(u_{j} -
u)|^{2}dt\leq C_{6}\int_{0}^{T}|u_{j} - u|^{2}\,dt<\varepsilon.
\end{eqnarray*}
This last one result, is consequence of the convergence
\eqref{e335}.\\
\\
These convergence implies that we can take limits in the approximate
problem \eqref{e311}-\eqref{e315}, and then to verify the conditions
$(i),$ $(ii),$ $(iii)$ and $(iv)$ of the Theorem.\\
\\
Now, we will make verify of the initial data and we prove the
uniqueness of solutions. Using the result of regularity we have that
\begin{eqnarray}
\label{e401}u,\,v\in C^{0}(0,\,T:\,L^{2}(\Omega))
\end{eqnarray}
In this form, makes sense we calculate $u(0)$ e $v(0).$ Let us
consider $\theta\in C^{1}(0,\,T:\,\mathbb{R}),$ with $\theta(0)=1$
and $\theta(T)=0.$ Since the convergence \eqref{e329} we have
\begin{eqnarray}
\label{e402}\int_{0}^{T}(u_{j}',\,z)\,\theta\,dt\longrightarrow
\int_{0}^{T}(u',\,z)\,\theta\,dt,\qquad z\in L^{2}(\Omega).
\end{eqnarray}
Performing integration by parts in \eqref{e402} we have
\begin{eqnarray}
\label{e403}-\;(u_{j}(0),\,z) - \int_{0}^{T}(u_{j},\,z)\,\theta'\,dt
\longrightarrow -\;(u(0),\,z) - \int_{0}^{T}(u,\,z)\,\theta'\,dt.
\end{eqnarray}
Using the convergence \eqref{e329} in \eqref{e403} we have
$(u_{j}(0),\,z)\longrightarrow (u(0),\,z),$ for all $z\in
H_{0}^{1}(\Omega).$ But $u_{j}(0)$ converges strong for $u_{0}$ in
$L^{2}(\Omega),$ consequently weak in $L^{2}(\Omega).$ Therefore
$(u_{j}(0),\,z)\longrightarrow (u_{0},\,z),$ for all $z\in
H_{0}^{1}(\Omega).$ From uniqueness of the limit we have
$(u(0),\,z)\longrightarrow(u_{0},\,z),$ for all $z\in
H_{0}^{1}(\Omega).$
Thus, $u(0)=u_{0}.$ In a similar way we can prove that $v(0)=v_{0}.$\\
\\
To finish this section we will show the uniqueness of solution. \\
\\
{\bf Theorem 3.2}(Uniqueness). {\it Let $(u_{0},\,v_{0})\in
L^{2}(\Omega)\times L^{2}(\Omega)$ and $0<T<+\infty,$ where the time
$T$ depends only $|u_{0}|_{L^{2}(\Omega)}$ and
$|v_{0}|_{L^{2}(\Omega)}.$ If \eqref{e206}-\eqref{e209} holds, then
there is at most one solution of \eqref{e301}-\eqref{e305} in
$L^{2}(0,\,T:\,H_{0}^{1}(\Omega))\cap
C([0,\,T):\,L^{2}(\Omega))\times
L^{2}(0,\,T:\,H_{0}^{1}(\Omega))\cap C([0,\,T):\,L^{2}(\Omega))$
with initial data
$(u(x,\,0),\,v(x,\,0))=(u_{0},\,v_{0}).$}\\
\\
{\it Proof.} Assume that $(u_{1},\,v_{1}),$ $(u_{2},\,v_{2})$ in
$L^{2}(0,\,T:\,H_{0}^{1}(\Omega))\cap
C([0,\,T):\,L^{2}(\Omega))\times
L^{2}(0,\,T:\,H_{0}^{1}(\Omega))\cap C([0,\,T):\,L^{2}(\Omega))$ are
two solutions of \eqref{e301}-\eqref{e305} with $u_{t},$ $v_{t}$ in
$L^{2}(0,\,T:\,H^{-1}(\Omega))\times L^{2}(0,\,T:\,H^{-1}(\Omega)),$
so all integrations below are justified and with the same initial
data, in fact $(u_{1} - u_{2})(x,\,0)\equiv 0$ and $(u_{2} -
u_{2})(x,\,0)\equiv 0.$ Then
\begin{eqnarray*}
\frac{d}{dt}u_{1} - a(l(u_{1}))\,\triangle u_{1} & = &
-f(u_{1},\,v_{1}) + g(u_{1},\,v_{1})\\
\frac{d}{dt}v_{1} - a(l(v_{1}))\,\triangle v_{1} & = &
f(u_{1},\,v_{1}) - g(u_{1},\,v_{1})
\end{eqnarray*}
and
\begin{eqnarray*}
\frac{d}{dt}u_{2} - a(l(u_{2}))\,\triangle u_{2} & = & -
f(u_{2},\,v_{2}) + g(u_{2},\,v_{2})\\
\frac{d}{dt}v_{2} - a(l(v_{2}))\,\triangle v_{2} & = &
f(u_{2},\,v_{2}) - g(u_{2},\,v_{2}).
\end{eqnarray*}
Then
\begin{eqnarray*}
\lefteqn{\frac{d}{dt}(u_{1} - u_{2},\,h_{1}) +
a(l(u_{1}))\int_{\Omega}\nabla u_{1}\cdot\nabla h_{1}\,dx -
a(l(u_{2}))\int_{\Omega}\nabla u_{2}\cdot\nabla h_{1}\,dx }
\\
& = & -(f(u_{1},\,v_{1}) - f(u_{2},\,v_{2}),\,h_{1}) +
(g(u_{1},\,v_{1}) -
g(u_{2},\,v_{2}),\,h_{1}),\qquad\forall\;h_{1}\in H_{0}^{1}(\Omega)
\end{eqnarray*}
\begin{eqnarray*}
\lefteqn{\frac{d}{dt}(v_{1} - v_{2},\,h_{2}) +
a(l(v_{1}))\int_{\Omega}\nabla v_{1}\cdot\nabla h_{2}\,dx -
a(l(v_{2})) \int_{\Omega}\nabla v_{2}\cdot\nabla h_{2}\,dx
}\\
& = & (f(u_{1},\,v_{1}) - f(u_{2},\,v_{2}),\,h_{2}) -
(g(u_{1},\,v_{1}) -
g(u_{2},\,v_{2}),\,h_{2}),\qquad\forall\;h_{2}\in H_{0}^{1}(\Omega).
\end{eqnarray*}
Using that $g(u_{1},\,v_{1}) - g(u_{2},\,v_{2})= \alpha\,(u_{1} -
u_{2}) - \alpha\,(v_{1} - v_{2})$ and that $f(u,\,v)=f(u - v)$ we
have
\begin{eqnarray*}
\lefteqn{\frac{d}{dt}(u_{1} - u_{2},\,h_{1}) +
a(l(u_{1}))\int_{\Omega}\nabla u_{1}\cdot\nabla h_{1}\,dx -
a(l(u_{2}))\int_{\Omega}\nabla u_{2}\cdot\nabla h_{1}\,dx }
\\
& = & -(f(u_{1} - v_{1}) - f(u_{2} - v_{2}),\,h_{1}) +
\alpha\,(u_{1} - u_{2},\,h_{1}) - \alpha\,(v_{1} -
v_{2},\,h_{1}),\quad\forall\;h_{1}\in H_{0}^{1}(\Omega)
\end{eqnarray*}
\begin{eqnarray*}
\lefteqn{\frac{d}{dt}(v_{1} - v_{2},\,h_{2}) +
a(l(v_{1}))\int_{\Omega}\nabla v_{1}\cdot\nabla h_{2}\,dx -
a(l(v_{2}))\int_{\Omega}\nabla v_{2}\cdot\nabla h_{2}\,dx
}\\
& = & (f(u_{1} - v_{1}) - f(u_{2} - v_{2}),\,h_{2}) - \alpha\,(u_{1}
- u_{2},\,h_{2}) + \alpha\,(v_{1} -
v_{2},\,h_{2}),\quad\forall\;h_{2}\in H_{0}^{1}(\Omega).
\end{eqnarray*}
On the other hand, let $h_{1}=u_{1} - u_{2}$ and $h_{2}=v_{1} -
v_{2}$ we obtain
\begin{eqnarray*}
\lefteqn{\frac{d}{dt}|u_{1} - u_{2}|^{2} +
a(l(u_{1}))\int_{\Omega}\nabla u_{1}\cdot\nabla (u_{1} - u_{2})\,dx
- a(l(u_{2}))\int_{\Omega}\nabla u_{2}\cdot\nabla (u_{1} -
u_{2})\,dx }
\nonumber \\
& = & -(f(u_{1} - v_{1}) - f(u_{2} - v_{2}),\,u_{1} - u_{2}) +
\alpha\,|u_{1} - u_{2}|^{2} - \alpha\,(v_{1} - v_{2},\,u_{1} -
u_{2})
\end{eqnarray*}
and
\begin{eqnarray*}
\lefteqn{\frac{d}{dt}|v_{1} - v_{2}|^{2} +
a(l(v_{1}))\int_{\Omega}\nabla v_{1}\cdot\nabla (v_{1} - v_{2})\,dx
- a(l(v_{2}))\int_{\Omega}\nabla v_{2}\cdot\nabla (v_{1} -
v_{2})\,dx
}\nonumber \\
& = & (f(u_{1} - v_{1}) - f(u_{2} - v_{2}),\,v_{1} - v_{2}) -
\alpha\,(u_{1} - u_{2},\,v_{1} - v_{2}) + \alpha\,|v_{1} -
v_{2}|^{2}.
\end{eqnarray*}
Hence
\begin{eqnarray*}
\lefteqn{\frac{d}{dt}|u_{1} - u_{2}|^{2} +
a(l(u_{1}))\,||u_{1}||^{2} + a(l(u_{2}))\,||u_{2}||^{2} +
[\,a(l(u_{1})) - a(l(u_{2})\,]\int_{\Omega}\nabla u_{1}\cdot\nabla
u_{2}\,dx}
\nonumber \\
& \leq & |f(u_{1} - v_{1}) - f(u_{2} - v_{2})|\;|u_{1} - u_{2}| +
\alpha\,|u_{1} - u_{2}|^{2} + \alpha\,|v_{1} - v_{2}|\;|u_{1} -
u_{2}|
\end{eqnarray*}
and
\begin{eqnarray*}
\lefteqn{\frac{d}{dt}|v_{1} - v_{2}|^{2} +
a(l(v_{1}))\,||v_{1}||^{2} + a(l(v_{2}))\,||v_{2}||^{2} +
[\,a(l(v_{1})) - a(l(v_{2})\,]\int_{\Omega}\nabla v_{1}\cdot\nabla
v_{2}\,dx
}\nonumber \\
& \leq & |f(u_{1} - v_{1}) - f(u_{2} - v_{2})|\;|v_{1} - v_{2}| +
\alpha\,|u_{1} - u_{2}|\;|v_{1} - v_{2}| + \alpha\,|v_{1} -
v_{2}|^{2}.
\end{eqnarray*}
Using \eqref{e206}-\eqref{e209} and the Young inequality we have
\begin{eqnarray*}
\lefteqn{\frac{d}{dt}|u_{1} - u_{2}|^{2} + m\,||u_{1}||^{2} +
m\,||u_{2}||^{2} }
\nonumber \\
& \leq & M_{1}\,|\,l(u_{1}) - l(u_{2})\,|\;||u_{1}||\,||u_{2}|| +
M_{3}\,|(u_{1} - u_{2}) - (v_{1} - v_{2})|\;|u_{1} - u_{2}|  \\
&  & +\;\alpha\,|u_{1} - u_{2}|^{2} + \frac{\alpha}{2}\,|v_{1} -
v_{2}|^{2} + \frac{\alpha}{2}\,|u_{1} - u_{2}|^{2}  \\
& \leq & M_{1}\,C_{1}\,|u_{1} - u_{2}|\;||u_{1}||\;||u_{2}|| +
M_{3}\,|u_{1} - u_{2}|^{2} + M_{3}\,|u_{1} - u_{2}|\;|v_{1} - v_{2}|   \\
&  & +\;\frac{3}{2}\alpha\,|u_{1} -
u_{2}|^{2} + \frac{\alpha}{2}\,|v_{1} - v_{2}|^{2} \\
& \leq & \frac{m}{2}\,||u_{1}||^{2} +
\frac{M_{1}^{2}\,C_{1}^{2}}{2m}\,||u_{2}||^{2}\,|u_{1} - u_{2}|^{2}
+ M_{3}\,|u_{1} - u_{2}|^{2} +
\frac{M_{3}}{2}\,|u_{1} - u_{2}|^{2}  \\
&  & +\;\frac{M_{3}}{2}\,|v_{1} - v_{2}|^{2} +
\frac{3}{2}\alpha\,|u_{1} -
u_{2}|^{2} + \frac{\alpha}{2}\,|v_{1} - v_{2}|^{2} \\
& \leq & \frac{m}{2}\,||u_{1}||^{2} +
\frac{M_{1}^{2}\,C_{1}^{2}}{2m}\,||u_{2}||^{2}\,|u_{1} - u_{2}|^{2}
+ \frac{3}{2}\left(M_{3} + \alpha\right)|u_{1} - u_{2}|^{2} +
\frac{1}{2}\left(M_{3} + \alpha\right)|v_{1} - v_{2}|^{2}
\end{eqnarray*}
and
\begin{eqnarray*}
\lefteqn{\frac{d}{dt}|v_{1} - v_{2}|^{2} + m\,||v_{1}||^{2} +
m\,||v_{2}||^{2}
}\nonumber \\
& \leq & M_{2}\,|\,l(v_{1}) - l(v_{2})\,|\,||v_{1}||\,||v_{2}|| +
M_{3}\,|(u_{1} - u_{2}) - (v_{1} - v_{2})|\;|v_{1} - v_{2}| \\
&  & +\;\frac{\alpha}{2}\,|u_{1} - u_{2}|^{2} +
\frac{\alpha}{2}\,|v_{1} - v_{2}|^{2} + \alpha\,|v_{1} -
v_{2}|^{2}\nonumber \\
& \leq & M_{2}\,C_{2}\,|v_{1} - v_{2}|\,||v_{1}||\,||v_{2}|| +
M_{3}\,|v_{1} - v_{2}|^{2} + M_{3}\,|u_{1} - u_{2}|\;|v_{1} - v_{2}|  \nonumber \\
&  & +\;\frac{3}{2}\alpha\,|v_{1} -
v_{2}|^{2} + \frac{\alpha}{2}\,|u_{1} - u_{2}|^{2} \nonumber \\
& \leq &  \frac{m}{2}\,||v_{1}||^{2} +
\frac{M_{2}^{2}\,C_{2}^{2}}{2m}\,||v_{2}||^{2}\,|v_{1} - v_{2}|^{2}
+ M_{3}\,|v_{1} - v_{2}|^{2} + \frac{M_{3}}{2}\,|v_{1} - v_{2}|^{2}  \nonumber \\
&  & +\;\frac{M_{3}}{2}\,|u_{1} - u_{2}|^{2} +
\frac{3}{2}\alpha\,|v_{1} - v_{2}|^{2} + \frac{\alpha}{2}\,|u_{1} -
u_{2}|^{2} \\
& \leq & \frac{m}{2}\,||v_{1}||^{2} +
\frac{M_{2}^{2}\,C_{2}^{2}}{2m}\,||v_{2}||^{2}\,|v_{1} - v_{2}|^{2}
+ \frac{3}{2}\left(M_{3} + \alpha\right)|v_{1} - v_{2}|^{2} +
\frac{1}{2}\left(M_{3} + \alpha\right)|u_{1} - u_{2}|^{2}
\end{eqnarray*}
Then
\begin{eqnarray}
\lefteqn{\frac{d}{dt}|u_{1} - u_{2}|^{2} +
\frac{m}{2}\,||u_{1}||^{2} + m\,||u_{2}||^{2} }
\nonumber \\
\label{e404}& \leq &
\frac{M_{1}^{2}\,C_{1}^{2}}{2m}\,||u_{2}||^{2}\,|u_{1} - u_{2}|^{2}
+ \frac{3}{2}\left(M_{3} + \alpha\right)|u_{1} - u_{2}|^{2} +
\frac{1}{2}\left(M_{3} + \alpha\right)|v_{1} - v_{2}|^{2}
\end{eqnarray}
and
\begin{eqnarray}
\lefteqn{\frac{d}{dt}|v_{1} - v_{2}|^{2} +
\frac{m}{2}\,||v_{1}||^{2} + m\,||v_{2}||^{2}
}\nonumber \\
\label{e405}& \leq &
\frac{M_{2}^{2}\,C_{2}^{2}}{2m}\,||v_{2}||^{2}\,|v_{1} - v_{2}|^{2}
+ \frac{3}{2}\left(M_{3} + \alpha\right)|v_{1} - v_{2}|^{2} +
\frac{1}{2}\left(M_{3} + \alpha\right)|u_{1} - u_{2}|^{2}
\end{eqnarray}
Adding \eqref{e404} with \eqref{e405} we obtain
\begin{eqnarray*}
\lefteqn{\frac{d}{dt}(\,|u_{1} - u_{2}|^{2} + |v_{1} - v_{2}|^{2}\,)
+ \frac{m}{2}\,||u_{1}||^{2} + m\,||u_{2}||^{2} +
\frac{m}{2}\,||v_{1}||^{2} + m\,||v_{1}||^{2} }
\nonumber \\
& \leq & \frac{M_{1}^{2}\,C_{1}^{2}}{2m}\,||u_{2}||^{2}\,|u_{1} -
u_{2}|^{2} + \frac{M_{2}^{2}\,C_{2}^{2}}{2m}\,||v_{2}||^{2}\,|v_{1}
- v_{2}|^{2} \nonumber \\
&  & +\;2\left(M_{3} + \alpha\right)|u_{1} - u_{2}|^{2} +
2\left(M_{3} + \alpha\right)|v_{1} - v_{2}|^{2} \\
& = & \left[\frac{M_{1}^{2}\,C_{1}^{2}}{2m}\,||u_{2}||^{2} +
2\,(M_{3} + \alpha)\right]|u_{1} - u_{2}|^{2} +
\left[\frac{M_{2}^{2}\,C_{2}^{2}}{2m}\,||v_{2}||^{2} + 2\,(M_{3} +
\alpha)\right]|v_{1} - v_{2}|^{2}
\end{eqnarray*}
We define
\begin{eqnarray*}
\varphi(t)=\frac{M_{1}^{2}\,C_{1}^{2}}{2m}\;||u_{2}||^{2} +
2\left(M_{3} + \alpha\right),\qquad \xi(t)=
\frac{M_{2}^{2}\,C_{2}^{2}}{2m}\;||v_{2}||^{2} + 2\left(M_{3} +
\alpha\right).
\end{eqnarray*}
Thus,
\begin{eqnarray*}
\frac{d}{dt}\left[\,|u_{1} - u_{2}|^{2} + |v_{1} -
v_{2}|^{2}\,\right]\leq \varphi(t)\,|u_{1} - u_{2}|^{2} +
\xi(t)|v_{1} - v_{2}|^{2}.
\end{eqnarray*}
Let  $\mathcal{R}(t)=\sup\{\varphi(t),\,\xi(t)\},$ then
$\mathcal{R}>0$ and
\begin{eqnarray}
\label{e412}\frac{d}{dt}\left(\,|u_{1} - u_{2}|^{2} + |v_{1} -
v_{2}|^{2}\,\right)\leq \mathcal{R}(t)\,\left(\,|u_{1} - u_{2}|^{2}
+ |v_{1} - v_{2}|^{2}\,\right).
\end{eqnarray}
Integrating \eqref{e412} over $t\in [0,\,T]$ and using that
$u_{1}(0)=u_{2}(0)$ and $v_{1}(0)=v_{2}(0),$ we obtain
\begin{eqnarray*}
|u_{1} - u_{2}|^{2} + |v_{1} - v_{2}|^{2} \leq
\int_{0}^{t}\mathcal{R}(t)\,(\,|u_{1} - u_{2}|^{2} + |v_{1} -
v_{2}|^{2}\,)\,dx.
\end{eqnarray*}
Let $\rho(t)=|u_{1} - u_{2}|^{2} + |v_{1} - v_{2}|^{2},$ then
\begin{eqnarray}
\label{e413}\rho(t)\leq \int_{0}^{t}\mathcal{R}(s)\,\rho(s)\,ds.
\end{eqnarray}
Applying Gronwall's inequality, we obtain
\begin{eqnarray*}
\rho(t)\leq 0.
\end{eqnarray*}
Therefore  $\rho\equiv 0,$ i. e.,
\begin{eqnarray*}
|u_{1} - u_{2}|^{2} + |v_{1} - v_{2}|^{2} = 0.
\end{eqnarray*}
Using the regularity of the solutions, the uniqueness follows.
\renewcommand{\theequation}{\thesection.\arabic{equation}}
\setcounter{equation}{0}\section{Exponential stability}In this
section we show that the total energy \eqref{e211} associated to
system \eqref{e201}-\eqref{e205} decay exponentially to zero as $t$
tends
to infinity. In what follows we will prove our main result:\\
\\
{\bf Theorem 4.1.} {\it Let $(u,\,v)$ be a solution of system
\eqref{e101}-\eqref{e105} given by the theorem 3.1 and theorem 3.2.
We suppose that $m > 2\,c_{p}\,(M_{1} + \alpha)> 0$,
where $c_p$ corresponds to the constant of the Poincar\'e inequality..
Then there exist positives constants $C$ and $\eta,$ such that}
\begin{eqnarray}
\label{e401}E(t) \leq  C\,E(0)\,e^{-\eta\,t}.
\end{eqnarray}
{\it Proof.} Multiplying equation \eqref{e201} by $u$ and
integrating over $x\in\Omega$ we have
\begin{eqnarray*}
\frac{1}{2}\frac{d}{dt}\int_{\Omega}|u|^{2}\,dx +
a(\ell(u))\int_{\Omega}|\nabla u|^{2}\,dx = -\int_{\Omega}f(u -
v)\,u\,dx + \alpha\int_{\Omega}(u - v)\,u\,dx.
\end{eqnarray*}
Multiplying equation \eqref{e202} by $v$ and integrating over
$x\in\Omega$ we have
\begin{eqnarray*}
\frac{1}{2}\frac{d}{dt}\int_{\Omega}|v|^{2}\,dx +
a(\ell(v))\int_{\Omega}|\nabla v|^{2}\,dx = \int_{\Omega}f(u -
v)\,v\,dx - \alpha\int_{\Omega}(u - v)\,v\,dx.
\end{eqnarray*}
Adding the expressions above and using \eqref{e208} we have
\begin{eqnarray*}
\frac{1}{2}\frac{d}{dt}\int_{\Omega}\left[\,|u|^{2} +
|v|^{2}\,\right]dx + m\int_{\Omega}\left[\,|\nabla u|^{2} + |\nabla
v|^{2}\,\right]dx \leq \int_{\Omega}|f(u - v)|\;|u - v|\,dx +
\alpha\int_{\Omega}|u - v|^{2}\,dx.
\end{eqnarray*}
Using \eqref{e206} follows that
\begin{eqnarray*}
\frac{d}{dt}E(t)\leq - m\int_{\Omega}\left[\,|\nabla u|^{2} +
|\nabla v|^{2}\,\right]\,dx + (M_{1} + \alpha)\int_{\Omega}|u -
v|^{2}\,dx.
\end{eqnarray*}
Applying the Poincar\'{e} inequality and the Young inequality, we
obtain
\begin{eqnarray*}
\frac{d}{dt}E(t)\leq - \frac{2m}{c_{p}}\,E(t) + 4\,(M_{1} +
\alpha)\,E(t).
\end{eqnarray*}
Now choosing  $m > 2\,c_{p}\,(M_{1} + \alpha)> 0$ follows that there
exists $C>0$ such that
\begin{eqnarray*}
\frac{d}{dt}E(t)\leq -\,C\,E(t).
\end{eqnarray*}
From we concludes that there exists $\eta>0$ such that
\begin{eqnarray*}
E(t)\leq \,C\,E(0)\,e^{-\eta\,t}.
\end{eqnarray*}
The proof follows.

\renewcommand{\theequation}{\thesection.\arabic{equation}}
\setcounter{equation}{0}\section{Numerical Results}

In this section we consider
a particular case of the nonlocal reaction diffusion equations \eqref{e201}-\eqref{e205}
in one dimensional  space ($n=1$), $\Omega=(0,1)$, and $l:L^{2}(\Omega)\rightarrow \mathbb{R}$  defined 
by $l(u)=\displaystyle\int_0^1u(x)dx$.
Then we approximate the solution of the system using
implicit finite differences.
The numerical scheme reads as following
\begin{eqnarray}
&&\frac{u^{k+1}_i - u^k_i}{\delta x} - a\left({\displaystyle\sum_{j=1}^J\delta x u^k_j}\right)\frac{u^{k+1}_{i+1}-2u^{k+1}_{i}+u^{k+1}_{i-1}}{\delta x^2} + f(u^k_i-v^k_i)
=
\alpha(u^k_i-v^k_i),\label{e501}\\
&&\frac{v^{k+1}_i - v^k_i}{\delta x} - a\left({\displaystyle\sum_{j=1}^J\delta x v^k_j}\right)\frac{v^{k+1}_{i+1}-2v^{k+1}_{i}+v^{k+1}_{i-1}}{\delta x^2} - f(u^k_i-v^k_i)
=
\alpha(v^k_i-u^k_i),\label{e502}\\
&& u^{k+1}_0=u^{k+1}_J=v^{k+1}_0=v^{k+1}_J=0\label{e503}
\end{eqnarray}
for $i=1,\ldots,J-1$ and $k=0,\ldots,K$, where
$\delta t = T/K$, $\delta x=1/J$, $x_i=i\delta x$, $i=0,\ldots,J$ and 
$t_k=k\delta t$, $k=0,\ldots,K$. In \eqref{e501}-\eqref{e503}, $u^k_i$ denotes the approximation of $u(t_k,x_i)$.
In order to solve the system \eqref{e501}-\eqref{e503}, we consider the initial condition approximated by
\begin{eqnarray}
 u^0_i=u_0(x_i),\quad\hbox{ and }\quad v^0_i=v_0(x_i),\qquad\hbox{for } i=0,\ldots,J.
\label{e504}
\end{eqnarray}
For each $k=0,\ldots,K$, the scheme \eqref{e501}-\eqref{e503} is equivalent to a linear system
with a tridiagonal matrix of $\mathbb{R}^{(J-1)\times(J-1)}$ which is positive definited
and then there exists a unique solution of \eqref{e501}-\eqref{e504}.

\begin{figure}[ht]
\hspace{-.8cm}
\begin{tabular}{cc}
\epsfig{file=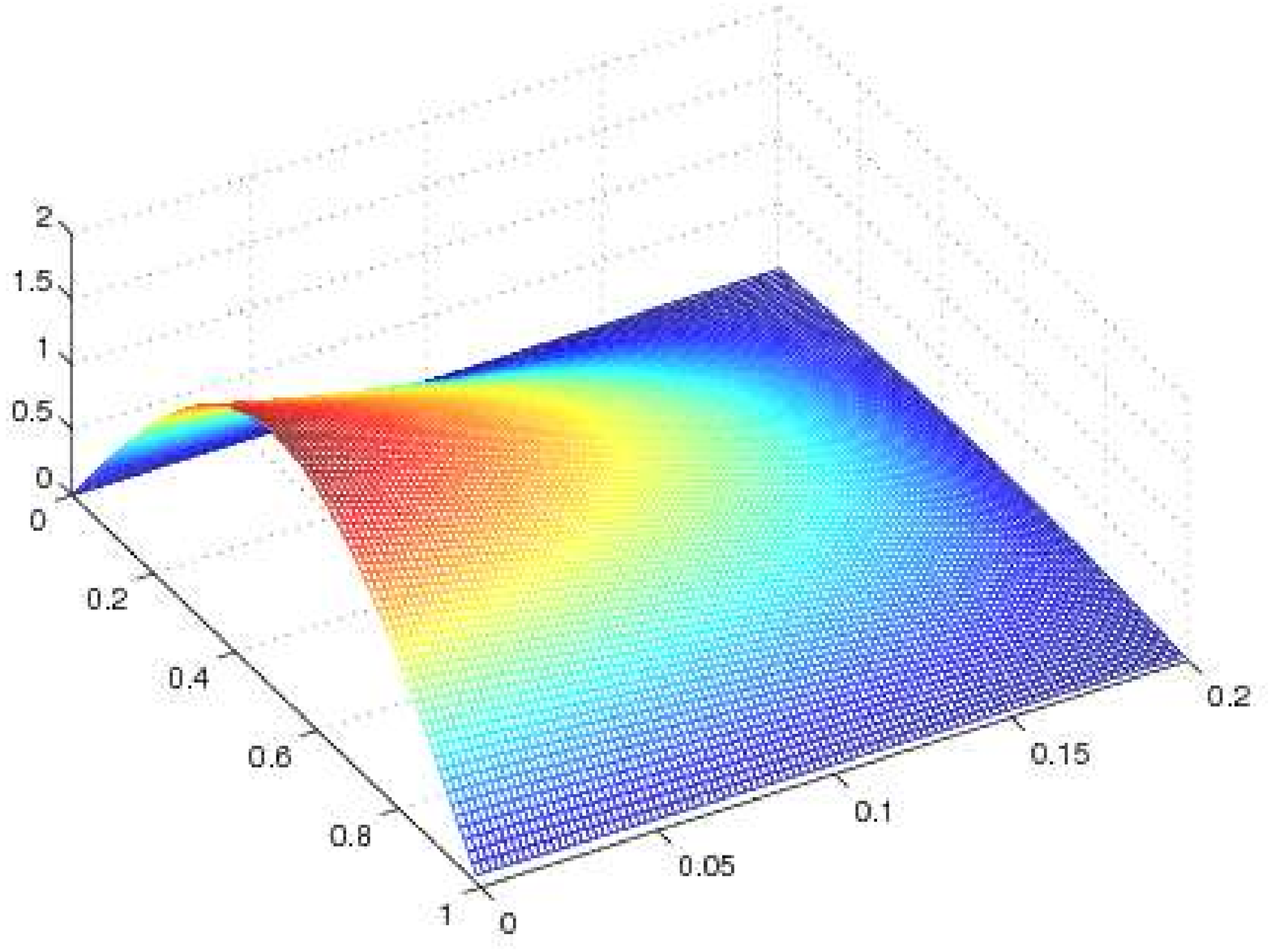,height=0.52\textwidth,width=0.52\textwidth} &
\epsfig{file=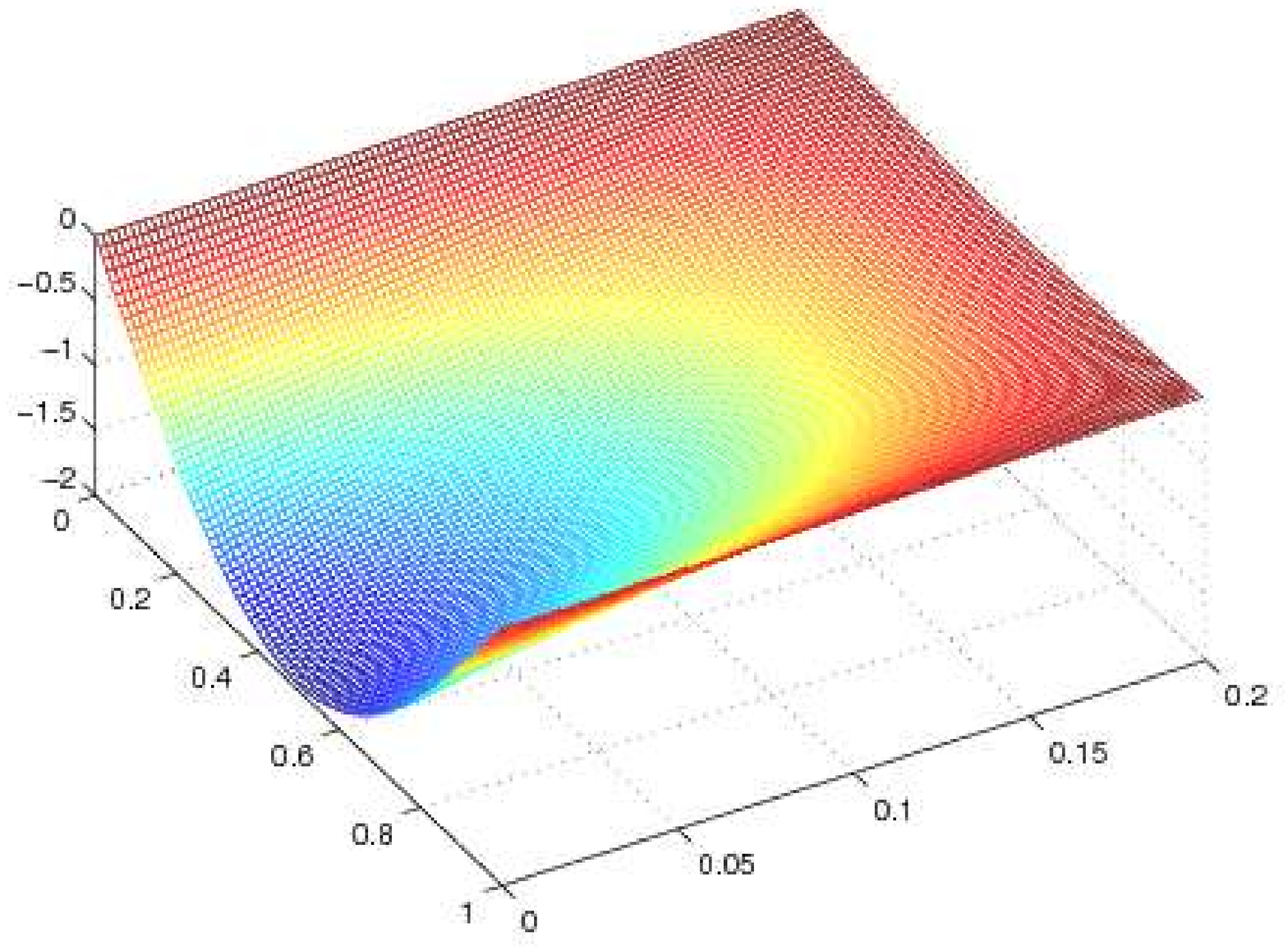,height=0.52\textwidth,width=0.52\textwidth}
\end{tabular}
\caption{Extinction of the population density $u(x,t)$ (left) and $v(x,t)$ (right) at time $T=0.2$.
}
\end{figure}

\begin{figure}[bt]
\hspace{-.8cm}
\begin{tabular}{cc}
\epsfig{file=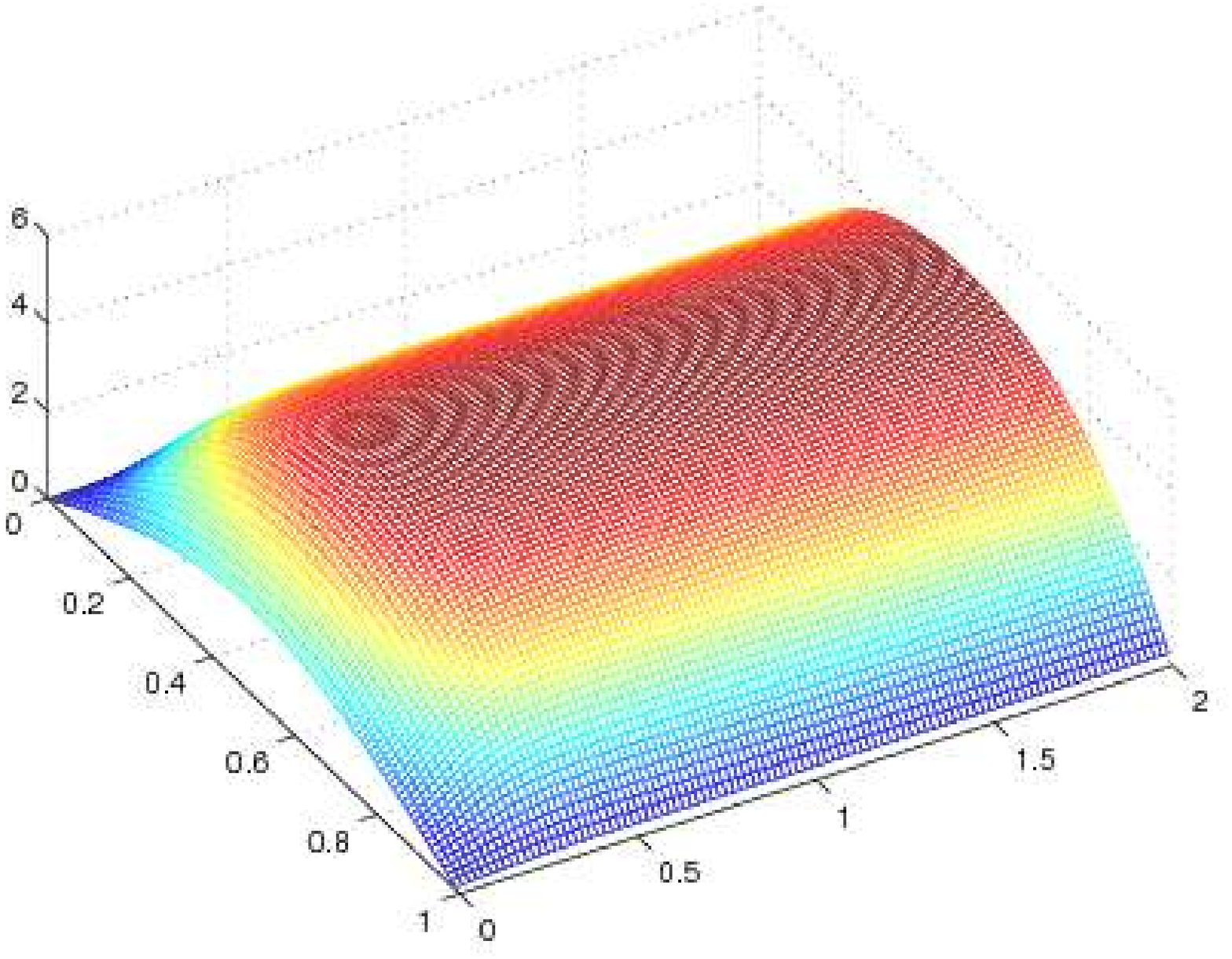,height=0.52\textwidth,width=0.52\textwidth} &
\epsfig{file=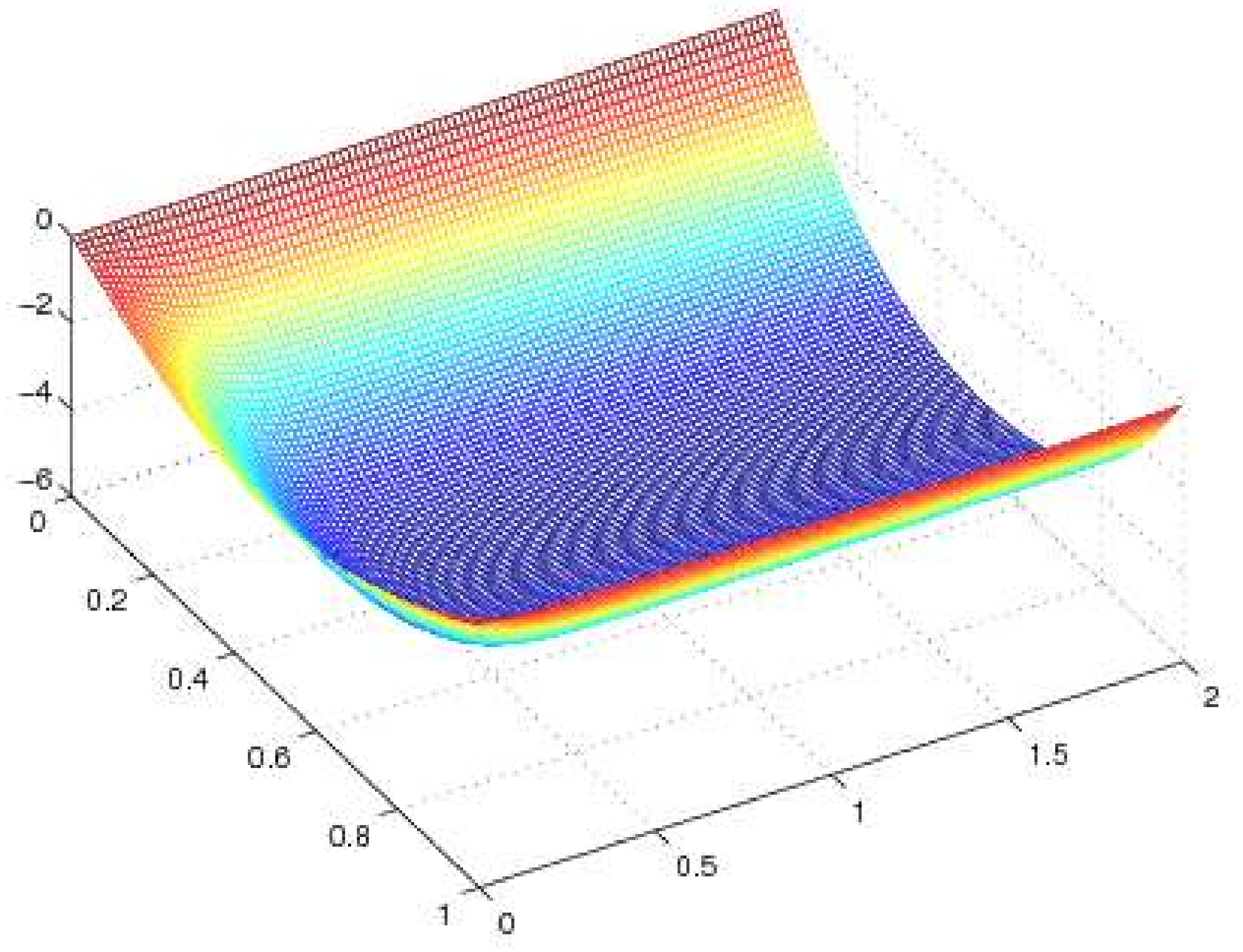,height=0.52\textwidth,width=0.52\textwidth}
\end{tabular}
\caption{Persistence of the population density $u(x,t)$ (left) and $v(x,t)$ (right).
}
\end{figure}

In order to compare the numerical behaviour of the solution with
theoretical and numerical behaviour of the solution for one equation (scalar case),
we take the same nonlinear reaction terms with similar parameters and initial conditions
of Ackle and Ke \cite{Azmy}. Nonlinear reaction and nonlocal diffusion are given by
\begin{eqnarray}
a(\xi) &:=& \max\left\{\varepsilon,\dfrac{1}{\left|\xi\right|}\right\} + m_0,\qquad \hbox{for all } \xi,
\label{e505}\\
f(w) - \alpha w &:=&  rw(\kappa-w),\qquad \hbox{for all } w,\label{e506}
\end{eqnarray}
where $\varepsilon$, $m_0$, $r$ and $\kappa$ are constant and positive parameters.
We remark that in the numerical example of Ackle and Ke \cite{Azmy}, the authors consider 
a nonlocal diffusion given by the expression $\dfrac{1}{\int_0^1u(x)dx}$.
In our case the parameter $\varepsilon$ is a very small parameter and
it plays a practical computational role to avoid the numerical overflow on the diffusion
when the extinction of the population occurs, that is when $u\approx 0$ or $v\approx 0$.
We consider a parameter $m_0$ to the numerical study of an extinsion case of population
and for a persistence case
of population. 
In fact, the exponential decaying of the energy \eqref{e401}, 
can be interpreted as the extinsion of two populations $u$ and $v$.
That is occurs when the hypothesis $m>c_p(M_1+\alpha)$ of the Theorem 4.1,
is verified. If $m$ is too small, the decaying of the energy is not guarantee
and a population persistence can be occur as we see in Figure 2. 

The initial condition is given by  $u_0(x)=\delta \sin(\pi x)$,
with $\delta=1.95$ (see \cite{Azmy}), and $v_0(x)= - u_0(x)=-\delta \sin(\pi x)$.
In this example, a negative $v_0$ has not a physical or biological sense, if $u$ and $v$ represent
densities of population, but we want to focus in the importance of the hypothesis  
of the Theorem 4.1, showing numerically that the exponential decaying of the energy does
no occur when the hypothesis is not verified.
We choose the parameter $\varepsilon=10^{-6}$ for the nonlocal
diffusion function \eqref{e505}, and we choose $r=1.0$ and $k=10.0$ for the reaction function \eqref{e506}.
The discretization is given by $J=10^4$ and $K=10^4$; we solve
the linear system \eqref{e501}-\eqref{e503} using  Thomas algorithm programming in Fortran90.

We consider here 2 simulations, the first one with $m_0=1.0$ (see Figure 1),
and the second one with $m_0=0.1$ (see Figure 2). The population persistence phenomenom
does not occur with a choice of a big amplitude $\delta$ for the initial condition as
it occurs in \cite{Azmy}. In fact we made simulations with different $\delta$ values and
the exponential decaying is always observed.

\subsection*{Acknowledgment}

MS  has been supported by  Fondecyt   project
\# 1070694 and
Fondap in Applied Mathematics (Project \# 15000001).
OVV acknowledge support by  Direcci\'on de Investigaci\'on de la Universidad del
B\'{\i}o-B\'{\i}o. DIPRODE projects, \# 061008 1/R.
Finally, OV and MS  acknowledge support by 
CNPq CONICYT Project, \# 490987/2005-2 (Brazil) and
\# 2005-075 (Chile).

\end{document}